\documentclass[11pt]{article}
\usepackage{amssymb,amsbsy,amsmath,amsfonts,amssymb,amscd,graphicx}
\usepackage[english,francais]{babel}
\newtheorem{theorem}{Theorem}[section]
\newtheorem{lemma}[theorem]{Lemma}
\newtheorem{proposition}[theorem]{Proposition}
\newtheorem{corollary}[theorem]{Corollary}

\newtheorem{remark}[theorem]{Remark}

\newcommand{\proof}{\noindent {\sl Proof.\/} }
\newcommand{\finproof}{\unskip\null\hfill$\square$\vskip 0.3cm}
\newcommand{\RR}{{\mathbb R}}
\newcommand{\NN}{{\mathbb N}}
\newcommand{\LL}{{L}}
\newcommand{\CC}{{\mathcal C}}
\newcommand{\FF}{{\mathcal E}}

\newcommand{\II}{{\mathcal I}}
\newcommand{\IL}{{\textsf I}}
\newcommand{\FL}{{\textsf E}}

\newcommand{\ee}{{\rm e}}
\newcommand{\dd}{\;{\rm d}}
\newcommand{\dx}{\;{\rm d}x}
\newcommand{\dy}{\;{\rm d}y}
\newcommand{\ds}{\;{\rm d}s}

\newcommand{\EE}{{\mathcal E}}
\newcommand{\EL}{{\textsf E}}
\newcommand{\DD}{{\mathcal D}}
\newcommand{\conv}{\!\circ\!}
\newcommand{\email}[1]{{\small E-mail: {\textsf {#1}}}}
\newcommand{\http}[1]{{\small Internet: {\textsf {#1}}}}

\begin{document}\selectlanguage{english}
\title{Large time asymptotics of the doubly nonlinear equation in the non-displacement convexity regime}

\author{
 Martial Agueh \footnote{ Department of Mathematics and Statistics
University of Victoria
PO Box 3060 STN CSC
Victoria, BC, Canada, V8W 3R4.
\email{agueh@math.uvic.ca},
\http{http://www.math.uvic.ca/$\sim$agueh/}},\kern8pt
Adrien Blanchet\footnote{GREMAQ, Universit\'e de Toulouse, Manufacture des tabacs, 21 all\'ee de Brienne, F-31000 Toulouse.
\email{adrien.blanchet@univ-tlse1.fr},
\http{http://www.ceremade.dauphine.fr/$\sim$blanchet/}}
\kern8pt\& \kern8pt Jos\'e A. Carrillo\footnote{ICREA (Instituci\'o Catalana de
Recerca i Estudis Avan\c cats) and Departament de
Mate\-m\`a\-tiques, Universitat Aut\`onoma de Barcelona, E-08193
Bellaterra, SPAIN. \email{carrillo@mat.uab.es},
\http{http://kinetic.mat.uab.es/$\sim$carrillo/}}}
\date{\today}
\maketitle
\begin{abstract}
We study the long-time asymptotics of the doubly nonlinear
diffusion equation
$\rho_t=\mbox{div}({|\nabla\rho^m|^{p-2}\nabla\rho^m})$ in
$\RR^n$, in the range $\frac{n-p}{n(p-1)}<m <\frac{n-p+1}{n(p-1)}$
and $1<p<\infty$ where the mass of the solution is conserved, but
the associated energy functional is not displacement convex. Using
a linearisation of the equation, we prove an $L^1$-algebraic decay
of the non-negative solution to a Barenblatt-type solution, and we
estimate its rate of convergence. We then derive the nonlinear
stability of the solution by means of some comparison method
between the nonlinear equation and its linearisation. Our results
cover the exponent interval $\frac{2n}{n+1}<p<\frac{2n+1}{n+1}$
where a rate of convergence towards self-similarity was still
unknown for the $p$-Laplacian equation.
\end{abstract}

\section{Introduction}
In this work, we consider the doubly nonlinear equation defined
for any $(t,x)\in (0,\infty)\times\RR^n$ by
\begin{equation}\label{dnl}
\left\{
\begin{array}{ll}
 \displaystyle \frac{\partial\rho}{\partial t}=\Delta_p \rho:=
 {\rm div} \left[\left|\nabla \rho^m \right|^{p-2}\nabla \rho^m \right]\,,
 &(x\in\RR^n, t>0)\vspace{.3cm}\\
\rho(t=0)=\rho_0\geq 0 \,,&(x\in\RR^n)
\end{array}
\right.
\end{equation}
with $1<p<\infty$, $0<m$ and $n\ge 3$. This class of equations
contains the linear diffusion equation, ($p=2, \;m=1$), commonly
known as the heat equation, $\partial_t\rho=\Delta \rho$ ; the
nonlinear diffusion equation $\partial_t\rho=\Delta \rho^m$, known
as the porous medium equation ($p=2,\; m>1$), or the fast
diffusion equation ($p=2,\; m<1$), and the gradient-dependent
diffusion equation,
$\partial_t\rho=\mbox{div}(|\nabla\rho|^{p-2}\nabla\rho):=\Delta_p\rho$,
that is, the $p$-Laplacian equation, ($p\neq 2, \;m=1$). When
$p\neq 2$ and $m\neq 1$, Eq.~\eqref{dnl} is called the doubly
nonlinear diffusion equation, due to the fact that its diffusion
term depends non-linearly on both the unknown density $\rho$, and
its gradient $\nabla \rho$. Such gradient-dependent diffusion equations appear
in several models in non-Newtonian fluids \cite{La}, in glaciology
\cite{Hu,CDDSV}, and in turbulent flows in porous media \cite{Le}.
For more details on these models, we refer to the recent monograph
of V\'azquez \cite{Vaz}, and the references therein.

Assuming that the initial data is integrable, $\rho_0\in
L^1(\RR^n)$, it is known that \eqref{dnl} has a unique solution
$\rho\in C\left([0,\infty), L^1(\RR^n)\right)$, with $\rho(t)\in
C^{1,\alpha}(\RR^n)$ for some $\alpha\in(0,1)$, see for instance
\cite{DiBen,DiH1,DiH2,L}. Moreover, starting with a non-negative
initial data, $\rho_0\geq 0$, it is known that the solution
$\rho(t)$ remains non-negative at all times. Furthermore for $n
\ge 3$, there exists a critical exponent,
\begin{equation*}
m_c:=\frac{n-p}{n(p-1)},
\end{equation*}
such that if $m>m_c$, then the mass of the solution is conserved,
$\int_{\RR^n}\rho(t)\dx=\int_{\RR^n}\rho_0\dx$, while if $m<m_c$,
the solution vanishes in finite time, see \cite{DiBen,Vaz} and the
references therein. In particular, for the $p$-Laplacian equation,
this corresponds to the critical $p$-exponent,
\begin{equation*}
p_c:=\frac{2n}{n+1},
\end{equation*}
above which  the mass of the solution is conserved, while the
solution disappears in finite time if $p<p_c$. Therefore, up to
renormalising the mass of $\rho_0$ to unity, we can assume without
loss of generality that, under the condition $m>m_c$, the solution
$\rho(t)$ of \eqref{dnl} is a density in $\RR^n$, for
all times $t\geq 0$.

By similarity and scaling, it can be shown that, above the
critical exponent $m_c$,  Eq.\eqref{dnl} has a unique self-similar
solution $\rho_{D_*}$, whose initial value is the Dirac mass at
the origin, that is, the fundamental solution of Eq.\eqref{dnl}.
In fact, among all the radially symmetric solutions of
\eqref{dnl}, this solution is the most concentrated whose initial
data have the same mass as $\rho_0$. It is called the Barenblatt
solution \cite{Ba}, and it is precisely:
\begin{equation*}
\rho_{D_*}(t,x)=\frac{1}{t^{n/\delta_p}}u_{D_*}\left(\frac{x}{t^{1/\delta_p}}\right),
\end{equation*}
where
\begin{equation*}
\delta_p:=n(p-1)(m-m_c)>0,
\end{equation*}
and
\begin{equation*}
u_{D_*}(y)=
\left\{
\begin{array}{lcl}
\displaystyle \frac{1}{\sigma}\,\mbox{exp}\left(-\frac{|p-1|^2}{p}|y|^{p/(p-1)}\right) &\mbox{if}& \displaystyle m=\frac{1}{p-1}\\
\displaystyle\left(D_*-\frac{m(p-1)-1}{mp}|y|^{p/(p-1)}\right)_+^{\frac{p-1}{m(p-1)-1}} &\mbox{if}&\displaystyle m\neq \frac{1}{p-1},
\end{array}
\right.
\end{equation*}
with $\sigma$ and $D_*$ are uniquely determined by the mass
conservation: $\|u_{D_*}\|_{L^1(\RR^n)}=
\|\rho_{D_*}(t)\|_{L^1(\RR^n)}= \|\rho_0\|_{L^1(\RR^n)}$.

When $p=2$ and $m>1-{2}/{n}$, the existence and uniqueness of the
Barenblatt solution was proved by Friedmann and Kamin in
\cite{FK}. Moreover, they showed that the solution $\rho(t)$ of
the Cauchy problem converges to $\rho_{D_*}(t)$ w.r.t. the
$L^1(\RR^n)$-norm, as $t\to\infty$, with no rates. Rates of
convergence were computed by Carrillo and Toscani \cite{CT} if
$m>1$, independently by Del Pino and Dolbeault \cite{DD1}, and
Otto \cite{Ot} if $m\geq 1-{1}/{n}$. The rates found in this range
were generically optimal. In the range $1-{2}/{n}<m<1-{1}/{n}$,
there were studies of the linearised problem by Carrillo,
Lederman, Markowich and Toscani \cite{CLMT}, and Denzler and McCann
\cite{DM}. These linearisations were useful to obtain rates of
decay for the nonlinear fast diffusion equation by Carrillo and
V\'azquez \cite{CV} and later by McCann and Slep{\v{c}}ev
\cite{MS}, and Kim and McCann \cite{KM}. The decay rates obtained
by using the linearisations are in general non optimal and is optimal in some sub-range, see~\cite{KM}.

When $p\neq 2$ and $m=1$, Kamin and V\'azquez \cite{KV} proved
existence and uniqueness of the Barenblatt solution
$\rho_{D_*}$ for the $p$-Laplacian equation when $p>p_c$,
along with an $L^1$-convergence of the solution $\rho(t)$ of the
Cauchy problem to $\rho_{D_*}(t)$, with no rates. Their proof
extends to the doubly nonlinear equation as long as $m>m_c$, see
\cite{Vaz}. Rates of convergence were computed by Del Pino and
Dolbeault \cite{DD2} when $p_c+{1}/{(n+1)}\leq p<n$ for the
$p$-Laplacian equation, but their rates are not optimal; see also
a similar result for the doubly nonlinear equation in \cite{DD3}.
In \cite{Ag-cras,Ag-na08}, Agueh generalises previous results by
deriving optimal rates for the convergence of the solution of the
Cauchy problem \eqref{dnl} to $\rho_{D_*}(t)$, for all $m\geq
m_c+{1}/{(n(p-1))}={(n-p+1)}/{(n(p-1))}$ and $p>1$. For instance,
when $p=2$, this condition coincides with the case $m\geq
1-{1}/{n}$, while for the $p$-Laplacian equation ($p\neq
2,\;m=1$), it corresponds to $p\geq
p_c+{1}/{(n+1)}={(2n+1)}/{(n+1)}$,  and therefore covers the range
$p\geq n$ left in \cite{DD2}, but not the remaining exponent
interval ${2n}/{(n+1)}<p<{(2n+1)}/{(n+1)}$. Similarly, for the
doubly nonlinear diffusion equation, the rate of convergence
remains unknown in the range
\begin{equation}\label{interval}
m_c< m<m_c+\frac{1}{n(p-1)}=\frac{n-p+1}{n(p-1)}.
\end{equation}
Indeed, the proof of \cite{Ag-na08} is based on optimal
transportation inequalities, which follows from the displacement
convexity \cite{Mc} of the energy functional associated with
\eqref{dnl}, that is, $H^F(\rho)=\int_{\RR^n} F[\rho]\dd x$, where

\begin{equation*}
F(x)=
\left\{
\begin{array}{lcl}
\displaystyle \frac{1}{p-1}x\ln x &\mbox{if}&\displaystyle  m=\frac{1}{p-1}\\ \\
\displaystyle \frac{mx^\gamma}{\gamma(\gamma -1)}, \;\gamma=m+\frac{p-2}{p-1} &\mbox{if}&\displaystyle  m\neq \frac{1}{p-1}.
\end{array}
\right.
\end{equation*}
This energy functional is displacement convex if and only if
$\gamma\geq 1-\frac{1}{n}$, or equivalent $m\geq
{(n-p+1)}/{(n(p-1))}$. This explains why the method of
\cite{Ag-na08} does not extend to the interval~\eqref{interval}.

The goal of this work is then precisely to derive a rate of
convergence w.r.t the $L^1(\RR^n)$-norm, of the non-negative
solution $\rho$ of the Cauchy problem~\eqref{dnl}, to the
Barenblatt solution $\rho_{D_*}(t)$, as $t\to\infty$, provided that
$m$ belongs to the remaining exponent interval~\eqref{interval},
that is,
\begin{equation}\label{range}
\frac{n-p}{n(p-1)}<m<\frac{n-p+1}{n(p-1)}.
\end{equation}
For convenience we rewrite the Cauchy problem~\eqref{dnl} as:
\begin{equation}\label{dnlfp}
\left\{
  \begin{array}{ll}
 \displaystyle    \frac{\partial\rho}{\partial t} =\mbox{div}\Big\{\rho\nabla c^*\left[\nabla\left(F'\conv \rho\right)\right]\Big\}, & (x\in\RR^n, \; t>0)\vspace{.3cm}\vspace{.3cm}\\
\rho(t=0) = \rho_0, & (x\in \RR^n),
  \end{array}
 \right.
\end{equation}
where $c^*(x)={|x|^p}/{p}$ is the Legendre transform of the convex function
\begin{equation*}
c(x)=\frac{|x|^q}{q}, \quad \frac{1}{p}+\frac{1}{q}=1.
\end{equation*}
By rescaling in time and space $\rho$ as follows:
\begin{equation}
\label{rescaledrho}
\rho(t,x)=\frac{1}{R(t)^n}u\left(\tau, y\right),
\end{equation}
where
\begin{equation}
\label{defR}
\tau=\ln R(t), \quad y=\frac{x}{R(t)},\quad  R(t)=(1+\delta_p t)^{1/\delta_p}, \quad \delta_p=(p-1)(nm+1)+1-n,
\end{equation}
it is easy to show that $\rho$ solves~\eqref{dnlfp} if and only if
$u$ solves the rescaled convection-diffusion equation
\begin{equation}\label{rescaledeq}
\left\{
  \begin{array}{ll}
 \displaystyle \frac{\partial u}{\partial \tau} =\mbox{div}\Big\{u\nabla c^*\left[\nabla\left(F'\conv u\right)\right]+uy\Big\} & (y\in \RR^n, \; \tau>0)\vspace{.3cm}\\
u(\tau=0) = \rho_0 & (y\in\RR^n) .
\end{array}
\right.
\end{equation}
Moreover, by conservation of mass there exists a unique $D_*$ such that the Barenblatt profile $u_{D_*}$ is the equilibrium
solution of~\eqref{rescaledeq}. Remark that in the considered
range of exponents, $m(p-1)-1<0$ and the Barenblatt profile is
simply given by
\begin{equation}\label{barenblattprofilefast}
u_{D_*}(y)=\left(D_*+\frac{1-\gamma}{m}c(y)\right)^{\frac{1}{\gamma-1}}.
\end{equation}
In fact, $u_{D_*}$ is the unique density function of same mass as $u_0$ which satisfies
on its support,
\begin{equation}\label{equilibriumeq}
\nabla\left(F' \conv u_{D_*}+c\right)=0.
\end{equation}

The main result of our paper is the following:

\begin{theorem}[Rates of convergence]\label{maintheo}
Let $m$ be in the range~\eqref{range} and $u_0$ a density
such that there exist positive constants $D_0> D_1$ for which
\begin{equation*}
u_{D_0}(x)\le \rho_0(x)=u_0(x)\le u_{D_1}(x)\quad\forall\;x
\in\RR^n\, . \tag{H1}
\end{equation*}
Consider $u$ a solution to~\eqref{rescaledeq} with initial data
$u_0$, there exists a unique $D_*$ such that $u(\tau)$ converges to the Barenblatt profile $u_{D_*}$ in $L^1(\RR^n)$. Moreover, there exist {a time $\tau_0$ and} two positive constants $\lambda$ and
$M=M(m,n,p,u_0,\tau_0)$ such that, {for any time $\tau>\tau_0$}
\begin{equation}\label{recaledL1decay}
\|u(\tau)-u_{D_*}\|_{L^1(\RR^n)} \leq
M\,e^{-\frac{\lambda}{2}\tau}\,.
\end{equation}
\end{theorem}
As a consequence, {for a time large enough} the corresponding solution $\rho(t)$
of~\eqref{dnl} converges to the Barenblatt solution
$\rho_{D_*}(t)$, algebraically fast in the $L^1$-norm, at the rate
$\lambda/(2\delta_p)$: there exist {a time $t_0$ and} a constant $C=C(m,n,p,\rho_0,t_0)$ such that, {for any time $t>t_0$} 
\begin{equation}\label{L1decay}
\|\rho(t)-\rho_{D_*}(t)\|_{L^1(\RR^n)} \leq C\,t^{-\lambda/(2\delta_p)},
\end{equation}
where $\delta_p=(p-1)(nm+1)+1-n$.

The main tool is the following {\sl relative free energy} with
respect to the Barenblatt solution $u_{D_*}$ defined by
\begin{equation}\label{energy}
  \FF[u|u_{D_*}] := \int_{\RR^n} \left[ F \conv u(y) -F \conv u_{D_*}(y) -F' \conv u_{D_*}(y) (u(x)-u_{D_*}(y) ) \right] \dy\;
\end{equation}
for any given $u\in \LL^1_+(\RR^n)$. Its derivative along the flow
of~\eqref{dnlfp} is formally given by
\begin{equation*}
 -\frac{\dd}{\dd \tau}  \FF[u(\tau)|u_{D_*}]=\II[u(\tau)|u_{D_*}]
\end{equation*}
where
\begin{equation*}
\II[u(\tau)|u_{D_*}] := \int_{\RR^n}\! u(\tau,y)\,\nabla\left( F' \conv u(\tau,y) +
 c(y)\right)\cdot\left( \nabla c^* \conv \nabla F' \conv u(\tau,y)+y\right)\dy\;.
\end{equation*}
In this paper, we prove that the relative entropy decays exponentially fast in the
form
\begin{equation}
  \label{energydecay}
  \FF[u(\tau)|u_{D_*}] \le \ee^{-\beta\,\tau}\FF[u_0|u_{D_*}],
\end{equation}
for some $\beta>0$. This is obtained in two steps. First, we
linearise~\eqref{rescaledeq} at the equilibrium solution $u_{D_*}$
by using the linear perturbation $u(\tau)=u_{D_*}+\epsilon
v(\tau)$, and we show that the linearised version of the relative
energy converges to $0$ exponentially fast, as in~\cite{CLMT}. For
that, we use the Hardy-Poincar\'e inequality recently established
by Blanchet, Bonforte, Dolbeault, Grillo and V\'azquez
in~\cite{BBDGV-cras}. Next, following the strategy
in~\cite{BBDGV-arch}, we try to compare the relative energy and
the dissipation of the relative energy --that is, the Fisher
information-- for both linearised and nonlinear equations, to
deduce the exponential decay~\eqref{energydecay} for the nonlinear
equation. The main differences with respect to \cite{BBDGV-arch}
lie in the fact that a direct relation between the linearised and
the nonlinear Fisher information is not clear due to the singular
characters at the origin of the weights when $1<p<2$. Therefore,
we are forced to use a sort of regularised linearised Fisher
information instead. Moreover, the control of the additional terms
appearing in the regularised entropy dissipation of the linearised
problem and in the relation between the entropy dissipations is
more involved in our case.

We note that, based on our computations (see
Remark~\ref{bakry-emery}), the Bakry-Emery approach used
in~\cite{CLMT}, which consists of differentiation twice the
relative energy $\EL[v(\tau)]$ to estimate the spectral gap at the
eigenvalue $0$, does not yield a positive result for our equation
when $1<p< 2$, and thus, a similar procedure to \cite{CV} for the
doubly nonlinear equations is not feasible. Moreover, the
Hardy-Poincar\'e inequality used here to establish the linear
stability is actually valid on a larger interval,
$m_*<m<m_c+\frac{1}{n(p-1)}$, which includes our interval
$m_c<m<m_c+\frac{1}{n(p-1)}$, as $m_*<m_c$, where 
$m_*:=\frac{n-2q}{n-q}+\frac{2-p}{p-1}$. Therefore, our
linearisation result extends naturally to the interval $m_*<m\leq
m_c$ where mass conservation for the nonlinear equation fails. In
this range, one needs to carefully define the right class of
initial data and a substitute of the Barenblatt solution, as done
in~\cite{BBDGV-arch} when $p=2$. Here, we will not follow this
path and we will restrict ourselves to the case
$m_c<m<m_c+\frac{1}{n(p-1)}$ where mass is conserved to
concentrate in the main new difficulties.

The paper is organised as follows. In
Section~\ref{sec:withoutrate}, we review and introduce the main
estimates on the solutions needed in the rest of the work. In
particular, we prove the convergence of the solution $u(\tau)$ of
\eqref{rescaledeq} to the equilibrium solution $u_{D_*}$ in
$C^1(\RR^n)$, as $\tau\to\infty$, with no rate. Then in
Section~\ref{sec:withratel}, we analyse a suitable linearised
problem for which we apply an entropy-entropy dissipation argument
based on Hardy-Poincar\'e inequalities. Finally,
Section~\ref{sec:withrate} is devoted to establish the exponential
decay of $u(\tau)$ to $u_{D_*}$ by the comparison between linear
and nonlinear relative entropy dissipations.
\section{Convergence without rate}\label{sec:withoutrate}
\setcounter{equation}{0}

Let us start by reviewing some well-known facts about the global
unique weak solutions associated to \eqref{dnl} in the range of
exponents considered. They conserve mass for all times, {\it i.e.},
\begin{equation*}
\int_{\RR^n} \rho(t,x) \dd x=\int_{\RR^n} \rho_0(x) \dd x\quad
\forall\;t\in(0,\infty)\;.
\end{equation*}
From now on, $D_*$ is the unique positive real such that
\begin{equation*}
  \int_{\RR^n} \rho_0(x)\dd x=\int_{\RR^n} u_{D_*}(x) \dd x\,.
\end{equation*}
Moreover, solutions of the Cauchy problem to \eqref{dnl} enjoy a
comparison principle and the $L^1$-contraction property. Due to
the change of variables \eqref{rescaledrho}, these properties hold
for the solution $u$ of the nonlinear Fokker-Planck
equation~\eqref{rescaledeq}. Since in the rest of this paper we
will only work with the scaled solutions of the nonlinear
Fokker-Planck equation \eqref{rescaledeq}, from now on we will use
$t$ instead of $\tau$ and $x$ instead of $y$ for the time and
position variables respectively. The quotient function
\begin{equation*}
  w(t,x):=\frac{u(t,x)}{u_{D_*}(x)}
\end{equation*}
is solution to
\begin{equation*}
  \frac{\partial w}{\partial t} = \frac{1}{u_{D_*}(x)} {\rm div}\left\{ w(t,x)\,u_{D_*}(x)
  \left[ \nabla c^* \conv \nabla F' \conv [w(t,x)\,u_{D_*}(x) ] -\nabla c^* \conv \nabla F' \conv u_{D_*}(x) \right]\right\}\;.
\end{equation*}
Define
\begin{equation*}
W_0:= \inf_{x\in\RR^n}\frac{u_{D_0}(x)}{u_{D_*}(x)}\le
\sup_{x\in\RR^n}\frac{u_{D_1}(x)}{u_{D_*}(x)}:=W_1\;.
\end{equation*}
A straightforward calculation gives
\begin{equation*}
W_0=\left(\frac{{D_*}}{D_0}\right)^{\frac{1}{1-\gamma}} \leq\quad
1\quad\leq\quad
\left(\frac{{D_*}}{D_1}\right)^{\frac{1}{1-\gamma}}:=W_1\;
\end{equation*}
with strict inequalities unless $\rho_0=u_{D_*}$. In terms of
$w_0=u_0/u_{D^*}$, the "sandwich" assumption on the initial data {\textsf
(H1)} of Theorem~\ref{maintheo} can be rewritten as follows: there exist positive constants
$D_0> D_1$ such that
\begin{equation*}
0<W_0\le \frac{u_{D_0}(x)}{u_{D_*}(x)}\le w_0(x)\le
\frac{u_{D_1}(x)}{u_{D_*}(x)}\leq W_1<\infty\quad\forall\;x
\in\RR^n\,.\tag{H1'}
\end{equation*}

\begin{remark}
Let us point out that the condition {\rm (H1)} or {\rm (H1')} in
the case of the fast diffusion equation ($p=2$) and in the corresponding
range, $1-2/n<m<1-1/n$, is not restrictive. In fact, as a
consequence of the Harnack inequalities proved in {\rm \cite{BV}},
the hypothesis {\rm (H1)} is satisfied by $\rho(t)$ for any $t>0$
with an initial data $u_0\in L^1_+(\RR^n)$. In the present case, a
similar Harnack inequality, not available in the literature, would
restrict the study of the asymptotic rates to this particular set
of initial data.
\end{remark}

As a consequence of the regularity theory of degenerate parabolic
equations \cite{DiBen}, we can control uniformly
$C^{1,\alpha}$-norms.

\begin{lemma}[Uniform $\CC^{1,\alpha}$-estimates]\label{ckesti}
Given a solution $u \in C([0,\infty);\LL^1(\RR^n))$
of~\eqref{rescaledeq} with initial data $u_0$ satisfying {\rm
(H1)}, then for any $t_0 \in (0,\infty)$,
\begin{equation*}
  \sup_{t \ge t_0}\|w(t)\|_{\CC^{1,\alpha}(\RR^n)} < \infty\;.
\end{equation*}
Moreover, there exists $C>0$ such that for any $x\in \RR^n \setminus B_1$
\begin{equation}\label{newbehinfty}
  \left|\nabla w (x)\right| \le C \frac{w(x)}{|x|}
\end{equation}
\end{lemma}
\proof Due to the comparison principle and the hypothesis (H1),
the function $u(t)$ is sandwiched between the two Barenblatt
profiles for all times, {\it i.e.},
$$
u_{D_0} \leq u(t) \leq u_{D_1} \qquad t \geq 0,
$$
and thus is uniformly bounded in $B_2$, the Euclidean ball of
radius 2, uniformly in $t\geq t_0>0$. Due to the regularity theory
of degenerate parabolic equations \cite{DiBen,DiH1,DiH2,L},
interior regularity estimates in the sense of $u(t)\in
\CC^{1,\alpha}(B_1)$, for any $0\leq\alpha<1$ hold uniformly in
$t\geq t_0>0$. Consider $w=u/u_{D_*}$, then $w$ is also bounded in $\CC^{1,\alpha}(B_1)$, for any $0\leq\alpha<1$ uniformly in $t\geq t_0>0$. To deal with large values of $x$ we
introduce, the rescaled function
\begin{equation*}
  u_\lambda(t,x):=\lambda^{p/(1-m)} u(t, \lambda x)
\end{equation*}
which is also solution to~\eqref{rescaledeq} but the annulus
$B_{2\lambda}/B_\lambda$ gets mapped into the annulus
$\Omega_1:=B_{2}/B_1$. Note that all derivatives of the rescaled
Barenblatt $u_{D_*/\lambda}$ are uniformly bounded from above and
below since
$$
D^\beta u_{D_*/\lambda}\to C \, D^\beta
\left(|x|^{\frac{q(p-1)}{m(p-1)-1}}\right) \qquad \mbox{uniformly in } \Omega_1
\mbox{ as } \lambda\to\infty
$$
for any multi-index $\beta\in \NN^d$. As a consequence, we get
that $w_\lambda(t)=u_\lambda(t)/u_{D_*/\lambda}$ is also uniformly
bounded from above and below in $\Omega_1$ uniformly in $\lambda
\ge 1$ and $t\ge t_0>0$. Again using the regularity theory of the
degenerate parabolic equation, we deduce
that the $\CC^{1,\alpha}$-norm of $w_\lambda(t)$ in $\Omega_1$ is
also uniformly bounded for $t\ge t_0$ and $\lambda \ge 1$ by a
constant $C_1$. Going back in the $\lambda$-scaling we find a
constant independent of $\lambda>1$ such that
\begin{equation}\label{jaifini}
  \frac{|\nabla w (t,\lambda x)|}{w(t,\lambda x)} \le
  \frac{C_1}{\lambda} \quad \mbox{in $(t_0,\infty)$}\;
\end{equation}
in $\Omega_1$, and thus, the $\CC^1$-norm of $w(t)$ in $\RR^n/B_1$
is uniformly bounded. Similar scaling argument applies to the
H\"older semi-norms. As a consequence of~\eqref{jaifini}
\begin{equation*}
  \left|\nabla w(t,\lambda x)\right| \le C_1 \frac{w(t,\lambda x)}{\lambda} \le 2\,C_1\frac{w(t,\lambda x)}{\lambda|x|}\;.
\end{equation*}
We thus obtain the desired result for any $y=\lambda x \in B_{2\lambda}\setminus B_{\lambda}$, and any $\lambda>0$.
\finproof

From this we can obtain the following result regarding the
evolution of the relative entropy to the stationary state.

\begin{proposition}[Entropy/entropy production]\label{energydissipation}
Let $u \in C([0,\infty);\LL^1(\RR^n))$ be a solution
of~\eqref{rescaledeq} for an initial data satisfying {\rm (H1)},
and consider the free energy $\EE$ defined by~\eqref{energy}. Its
derivative along the flow of~\eqref{rescaledeq} is:
\begin{equation*}
  \frac{\dd}{\dd t}\EE[u(t)|u_{D_*}]=-\II\left[u(t)|u_{D_*}\right]
\end{equation*}
where
\begin{equation*}
\II\left[u(t)|u_{D_*}\right]\!:=\!\!\int_{\RR^n}\!\!\!
u(t)\!\left[\nabla\left(F' \conv u(t)-F' \conv
u_{D_*}\right)\right]\!\cdot \!\left[\nabla c^*\left(\nabla
\left(F' \conv u(t)\right)\right) - \nabla c^*\left(\nabla
\left(F' \conv u_{D_*}\right)\right)\right]\!\!\!\dy
\end{equation*}
is the relative Fisher information of $u(t)$ w.r.t. $u_{D_*}$. Moreover,
$\II\left[u(t)|u_{D_*}\right]=0$ if and only if $u=u_{D_*}$.
\end{proposition}
\proof By performing formally integration by parts, we get
\begin{eqnarray*}
\frac{\dd}{\dd t}\EE[u(t)|u_{D_*}] &=& \int_{\RR^n} \left[F' \conv u(t)-F'\conv u_{D^*}\right]\mbox{div}
\Big\{u(t)\nabla c^*\left[\nabla\left(F' \conv u(t)\right)\right]+u(t)y\Big\}\dy \\
&=& -\int_{\RR^n} u(t)\nabla\left[F' \conv u(t)-F'\conv u_{D^*}\right]\cdot
\Big\{\nabla c^*\left[\nabla\left(F' \conv
u(t)\right)\right]+y\Big\}\dy.
\end{eqnarray*}
The above energy dissipation follows using that $u_{D_*}$
satisfies \eqref{equilibriumeq} and $\nabla c^* \conv \nabla
c={\rm id}$. This integration by parts can be justified using
Lemma \ref{ckesti} by a standard argument introducing a cut-off
function like in \cite[Proposition 2.6]{BBDGV-arch}. Since the
arguments are exactly equal, we do not perform any further
details. By the convexity of $c^*$,
\begin{equation}\label{concavityc}
 \left[ \nabla c^*(a)-\nabla c^*(b)\right] \cdot \left( a-b\right) \ge 0
\end{equation}
with equality if and only if $a=b$. So the Fisher information is
non-negative and zero if and only if $u$ and $u_{D_*}$ have the same mass and such that $\nabla\left(F' \conv
u(\tau)-F' \conv u_{D_*}\right)=0$, {\it i.e.} $u=u_{D_*}$.
\finproof

With these ingredients, we can obtain a first result of
convergence toward stationary states.

\begin{lemma}[Uniform convergence]\label{prop:cvrgcewithoutrate}
Let $u \in C([0,\infty);\LL^1(\RR^n))$ be a solution
of~\eqref{rescaledeq} for an initial data satisfying {\rm (H1)},
then $\lim_{t\to\infty}w(t,x)=1$ uniformly in compact sets of
$\RR^n$. \end{lemma}
\proof Define $u^h(t,x):=u(h+t,x)$, for any given $h>0$ and
$t\in[0,1]$. It is also well-known \cite{DiBen} that equi-bounded
set of solutions of \eqref{dnl} are equi-continuous in time. This
property carries over to $u(t)$ by the change of variables in
\eqref{rescaledrho}. This fact together with the uniform bounds in
$C^{1,\alpha}$ obtained in Lemma~\ref{ckesti} and the
Ascoli-Arzel\'a theorem implies that for any sequence
$(h_{n})_{n\in\NN}$ there exists a sub-sequence
$(h_{n})_{n\in\NN}$, denoted with the same index, such that
$\{u^{h_n}\}_{n\in\NN}$ converges to a function $u_{\infty}$
uniformly in compact sets of $[0,1]\times\RR^n$, and moreover,
$u_\infty(t)\in C^{1,\alpha}(\RR^n)$ for all $t\in [0,1]$. Since
$\FF[u(t)|u_{D_*}]$ is non-increasing and positive and
\[
\FF[u(h_{n})|u_{D_*}]-\FF[u(h_{n}+1))|u_{D_*}] =
\int_{h_{n}}^{h_{n}+1}\II[u(s)|u_{D_*}]\ds
=\int_{0}^{1}\II[u(s+h_{n})|u_{D_*}]\ds\;,
\]
the function $t \mapsto \II[u^{h_{n}}(t)|u_{D_*}]$ is integrable
on $[0,1]$ and converges to zero as $n\to\infty$.
By~\eqref{concavityc}, $\II$ is non-negative. By Fatou's lemma we
have
\begin{equation*}
 \int_{\RR^n}  \lim_{n \to \infty} u^{h_{n}}(t,x)\,\nabla
 \left( F' \conv u^{h_{n}}(t,x) + c(x)\right)\cdot
 \left( \nabla c^*\left[ \nabla F' \conv u^{h_{n}}(t,x)+
 \nabla c^* \conv \nabla c(x)\right]\right)\dx=0\;.
\end{equation*}
As a consequence of~\eqref{concavityc}, $u_{\infty}$ satisfies
$\nabla[F' \conv u_{\infty}(x) + c(x)] = 0$, from which
$u_\infty=u_{D}$ for some $D>0$. By conservation of mass $D=D_*$.
Since the limit of all the convergent sub-sequences is uniquely
determined by $u_{D_*}$, the result is proved. \finproof

\begin{proposition}[Convergence in $\LL^p$-spaces]
Let $u \in C([0,\infty);\LL^1(\RR^n))$ be a solution of the scaled
doubly-nonlinear equation~\eqref{rescaledeq} for an initial data
satisfying {\rm (H1)}, then $\lim_{t \to
\infty}\|u(t)-u_{D_*}\|_p=0$, for any $p \in [1,\infty]$.
\end{proposition}
\proof By Lemma~\ref{prop:cvrgcewithoutrate}, $\lim_{t \to
\infty}\left|u(t,x)-u_{D_*}(x)\right|=0$ for any $x\in\RR^n$.
Moreover, by assumptions~{\rm (H1)}, for $|x|$ large enough
\[
\left|u(t)-u_{D_*}\right|\leq\max\big\{\left|u_{D_0}-u_{D_*}\right|,\;
\left|u_{D_1}-u_{D_*}\right|\big\}=O\left(|x|^{-q(2-\gamma)/(1-\gamma)}\right)\;.
\]
So the difference between $|u(t)-u_{D_*}|^\theta$ is in
$\LL^1(\RR^n)$ if $\theta >  \Theta(p,m)$ with
$$
\Theta(p,m):=\frac{n(1-\gamma)}{q(2-\gamma)}.
$$
It is easy to check that $\Theta(p,m)$ is a decreasing function of $\gamma$ and so of
$m$. Since $q>2$, in the range of exponents \eqref{range}, we have
$$
\Theta(p,m) \le \Theta(p,m_c)=\frac{n}{n+q}<1 \,.
$$
By Lebesgue's dominated convergence theorem, it
implies that $u(t)$ converges to $u_{D_*}$ in $\LL^\theta(\RR^n)$,
for any $\theta\in [1,\infty)$. Finally, we use the following
interpolation lemma, due to Nirenberg,~\cite[p. 126]{N}:
$$
\|f\|_{\infty} \le C\,\|f\|_{C^1(\RR^n)}^{\frac{n}{n+2}} \;
\|f\|_{2}^{\frac{2}{n+2}}\quad\forall\;f\in C(\RR^n)\;,
$$
for $f=u(t)-u_{D_*}$ together with Lemma~\ref{ckesti} to obtain the
result in the uniform norm. \finproof

\begin{remark}
In contrast with {\rm\cite{BBDGV-arch}}, we do not generally have
the convergence in $C^{1,\alpha}(\RR^n)$.
\end{remark}

\section{Linear stability}\label{sec:withratel}
\setcounter{equation}{0}

To prove the decay \eqref{energydecay} of $u(t)$ to $u_{D_*}$ in
the energy form, it is sufficient to establish the following
logarithmic Sobolev type inequality:
\begin{equation}\label{logsob}
 \EE\left[u|u_{D^*}\right] \leq \frac{1}{\beta} \II\left[u|u_{D^*}\right],
\end{equation}
for some $\beta>0$ and $u\in C^{1,\alpha}(\RR^n)$ verifying (H1).
Indeed, \eqref{logsob} combined with
Proposition~\ref{energydissipation} yield
\begin{equation*}
  \frac{\dd}{\dd t}\EE\left[u(t)|u_{D^*}\right] \leq -\beta \EE\left[u(t)|u_{D^*}\right],
\end{equation*}
and this leads to the energy decay~\eqref{energydecay} by a simple
Gronwall argument. To prove \eqref{logsob}, we will first show a
linearised version of this inequality, by considering the
perturbation $u(t)=u_{D_*} +\epsilon v(t)$ of the solution $u(t)$
to \eqref{rescaledeq}. This will lead to the convergence of $v(t)$
to 0 in relative entropy for the linearised equation
of~\eqref{rescaledeq}, as we will show below. Next section will be
devoted to compare the relative entropy and Fisher information
in~\eqref{logsob} with their linearised analogues.\medskip

For clarity in our exposition, let us start by formally deriving the form
of the linearised logarithmic Sobolev inequality that we will be
dealing with below. Using the perturbation $u=u_{D^*} +\epsilon v$
and the second order Taylor expansion of $F\left(u_{D^*} +\epsilon
v\right)$ at $\epsilon=0$ on the expression of the relative
entropy \eqref{energy}, we have that
\begin{equation*}
  F \conv u-F \conv u_{D^*}=\epsilon v F'(u_{D^*}) + \frac{\epsilon^2}{2}v^2F^{\prime\prime}\conv u_{D^*} + O(\epsilon^3),
\end{equation*}
and then $\EE\left[u|u_{D^*}\right]$ linearises as:
\begin{equation*}
\EE\left[u|u_{D^*}\right] = \frac{\epsilon^2}{2} \int_{\RR^n}
v^2F^{\prime\prime} \conv u_{D^*} + O(\epsilon^3).
\end{equation*}
Let us hence introduce the weighted $L^2$-norm:
\begin{equation}\label{neweq}
\EL\left[v\right]=\frac{1}{2}\int_{\RR^n}
v^2(x)\,F^{\prime\prime}\conv u_{D^*}(x)\dd x,
\end{equation}
which will play the role of the linearised relative entropy.

Concerning the linearisation of the Fisher information, from the first order Taylor expansion of $F'
\conv u=F'\left(u_{D_*}+\epsilon v\right)$ at $\epsilon=0$, we
have
\begin{equation}\label{eq2sec3}
B:=\nabla\left[F' \conv u \right] =A+\epsilon W+O(\epsilon^2).
\end{equation}
with
\begin{equation*}
  A:=\nabla\left[ F' \conv u_{D^*}\right]=\nabla c \quad \mbox{and} \quad W:=\nabla\left[vF^{\prime\prime} \conv u_{D^*}\right].
\end{equation*}
Then using that $\nabla c^*(z)=z|z|^{p-2}$, we obtain
\begin{equation}\label{eq3sec3}
\nabla c^*(B) = \nabla c^*(A) + \epsilon |A|^{p-2} W +
\epsilon(p-2) |A|^{p-4} (A \cdot W) A +O(\epsilon^2).
\end{equation}
Combining~\eqref{eq2sec3} and~\eqref{eq3sec3}, we
see that $\II\left[u|u_{D^*}\right]$ formally linearises as:
\begin{equation*}
\II\left[u|u_{D^*}\right] = \epsilon^2 \int_{\RR^n} u_{D^*}
|A|^{p-2} |W|^2\dx + \epsilon^2(p-2)\int_{\RR^n} u_{D^*} |A|^{p-4}
(A\cdot W)^2\dx +O(\epsilon^3).
\end{equation*}
Hence, for $\epsilon$ small enough, the logarithmic Sobolev inequality (\ref{logsob}) linearises as
\begin{align}\label{eq5sec3}
\frac{\beta}{2}\int_{\RR^n} v^2 F^{\prime\prime} \conv u_{D^*}
\leq& \int_{\RR^n} u_{D^*} |\nabla c|^{p-2}
|\nabla\left[v F^{\prime\prime} \conv u_{D^*}\right] |^2\dx \\
 & +(p-2)\int_{\RR^n} u_{D^*} |\nabla c|^{p-4} \left[\nabla c
\cdot \nabla\left(v F^{\prime\prime} \conv
u_{D^*}\right)\right]^2\dx. \nonumber
\end{align}
It will be shown below that the l.h.s of~\eqref{eq5sec3} is a
Lyapunov function -- and the relative entropy -- for the
linearised equation of~\eqref{rescaledeq}, and the r.h.s
of~\eqref{eq5sec3} corresponds to the dissipation of this relative
entropy, up to a constant.\medskip

Let $u$ be the solution
of~\eqref{rescaledeq}, and consider the small perturbation
\begin{equation}\label{pertubation}
u(t)=u_{D^*} +\epsilon v(t)
\end{equation}
of $u_{D^*}$, where $\epsilon>0$ is small, and $v(t)\in
C^{1,\alpha}(\RR^n)$ for some $\alpha\in (0,1)$. Because of the
mass-conservation, we have that
\begin{equation*}
\int_{\RR^n} v(t,x)\dx=0, \quad \forall t\geq 0.
\end{equation*}
Moreover~\eqref{pertubation} implies that
\begin{equation}\label{eq1sec3.1}
\frac{\partial u}{\partial t} =\epsilon\frac{\partial v}{\partial
t}.
\end{equation}
On the other hand, using~\eqref{equilibriumeq}, we have that
$\nabla c^*\left(A\right)=\nabla c^*[-\nabla c(x)]=-x$, and
then~\eqref{eq3sec3} gives that
\begin{equation}\label{eq2sec3.1}
u\left[\nabla c^*\left(B\right) +x\right] = \epsilon u_{D^*}
\left[ |A|^{p-2} W + (p-2)|A|^{p-4} (A \cdot W) W \right]
+O(\epsilon^2).
\end{equation}
Inserting~\eqref{eq1sec3.1}-\eqref{eq2sec3.1}
into~\eqref{rescaledeq}, we formally obtain after simplifying by
$\epsilon$ and then setting $\epsilon=0$, that the linearised
problem to~\eqref{rescaledeq} is:
\begin{equation}\label{linearizedeq}
\left\{
  \begin{array}{ll}
\displaystyle\frac{\partial v}{\partial t} =\mbox{div}\Big\{u_{D^*}\left(|A|^{p-2}W+(p-2)|A|^{p-4}(A\cdot W)A \right)\Big\}
& (x\in \RR^n, \; t>0)\vspace{.3cm}\\
v(t=0) = v_0 & (x\in\RR^n),
\end{array}
\right.
\end{equation}
with $v_0\in L^1(\RR^n)$ of zero average. We can easily check that
Eq.~\eqref{linearizedeq} has the linearised relative entropy~\eqref{neweq} as
Lyapunov functional. Actually, differentiating
$\EL\left[v(t)\right]$ along a solution $v$
of~\eqref{linearizedeq}, we formally have by a straightforward computation,
that
\begin{equation}\label{linearizeddissipation}
\frac{\dd }{\dd t}\EL\left[v(t)\right] =
-\left(\IL\left[v(t)\right] + (p-2)\IL_0\left[v(t)\right]\right),
\end{equation}
where
\begin{equation*}
\IL\left[v(t)\right]=\int_{\RR^n} |W(t)|^2u_{D^*}|A|^{p-2}\dd x
\quad\mbox{and}\quad \IL_0\left[v(t)\right]=\int_{\RR^n}
\left(A\cdot W(t)\right)^2u_{D^*}|A|^{p-4} \dd x.
\end{equation*}
The Cauchy-Schwarz inequality implies that $|A|^{p-4}\left(A\cdot
W(t)\right) \leq  |A|^{p-2}|W(t)|^2$, and as a consequence,
$\IL_0\left[v(t)\right] \leq \IL\left[v(t)\right]$. Using $1<p<2$,
we have
\begin{equation}\label{I}
 \IL\left[v(t)\right] + (p-2)\IL_0\left[v(t)\right] \geq (p-1)\IL\left[v(t)\right]\geq 0.
\end{equation}
In case $p>2$, it is direct that $\IL\left[v(t)\right] +
(p-2)\IL_0\left[v(t)\right] \geq \IL\left[v(t)\right]\geq 0$. From
these estimates, the dissipation \eqref{linearizeddissipation} and
$|A(x)|=|\nabla\left[ F' \conv u_{D^*}\right](x)|>0$ for all
$x\in\RR^n$, we readily formally conclude that the unique steady
state is the zero solution.\medskip

The objective of the rest of this section is to show the following
asymptotic exponential relaxation of the linearised equation
\eqref{linearizedeq}:
\begin{theorem}\label{theolinearizeddecay}
Let $m$ satisfying~\eqref{range} and $v_0\in L^1(\RR^n)$ with zero
average. Consider $v(t)$ the solution to \eqref{linearizedeq} with initial data $v_0$. There exists a constant $\beta>0$ such that
\begin{equation}\label{linearizeddecay}
\EL\left[v(t)\right] \leq e^{-\beta t}\EL[v_0].
\end{equation}
\end{theorem}

Let us concentrate first in the case $1<p<2$. To derive this
exponential rate of convergence, we establish the following
functional inequality
\begin{equation*}
\EL\left[v\right] \leq \frac{1}{\beta} \left(\IL\left[v\right] +
(p-2)\IL_0\left[v\right]\right), \quad \beta>0,
\end{equation*}
for all $v\in C^{1,\alpha}(\RR^n)$ with zero average. This inequality corresponds to the formal linearisation~\eqref{eq5sec3} of the logarithmic Sobolev
inequality~\eqref{logsob}. In fact, because of~\eqref{I}, it is
sufficient to prove the following linearised logarithmic Sobolev
type inequality:
\begin{equation}\label{linearizedlogsob}
\EL\left[v\right] \leq \frac{p-1}{\beta}\IL\left[v\right], \quad
\beta>0
\end{equation}
for all $v\in C^{1,\alpha}(\RR^n)$ with zero average. This is
equivalent to show the Hardy-Poincar\'e type inequality:
\begin{equation}\label{eq5sec3.1}
\int_{\RR^n} g^2\dd\mu(x) \leq \tilde\beta\int_{\RR^n} |\nabla
g|^2\dd\nu(x),
\end{equation}
for some $\tilde\beta>0$, where $g$ is any function satisfying
$\int_{\RR^n} g\dd\mu(x)=0$, and
\begin{equation}\label{defmunu}
\dd\mu(x)=\frac{\dd x}{F^{\prime\prime}\conv u_{D^*}(x)}, \quad
\dd\nu(x)=u_{D^*}(x)|\nabla\left(F'\conv
u_{D^*}(x)\right)|^{p-2}\dd x.
\end{equation}
Note that the functions $v$ and $g$
in~\eqref{linearizedlogsob}-\eqref{eq5sec3.1} are related by
$g=v\,F^{\prime\prime}\conv u_{D^*}$, and
$\tilde\beta=2(p-1)/\beta$. The inequality \eqref{eq5sec3.1} is
also enough to prove the needed inequality in the case $p>2$ since
\eqref{linearizeddissipation} implies that
$$
\frac{\dd }{\dd t}\EL\left[v(t)\right] \leq -
\IL\left[v(t)\right].
$$

\begin{lemma}[Hardy-Poincar\'e type inequality]\label{lemmahardypoincare}
Let $m, n, p$ be such that $1<p<\infty$ and
$m_c<m<\frac{n-p+1}{n(p-1)}$. Then, there exits a constant
$\tilde\beta >0$ such that
\begin{equation*}
\int_{\RR^n} g^2\dd\mu(x) \leq \tilde\beta\int_{\RR^n} |\nabla
g|^2\dd\nu(x),
\end{equation*}
for any function $g\in C^{1,\tilde\alpha}(\RR^n)$ satisfying
$\int_{\RR^n}g\dd\mu(x)=0$ with $0<\tilde\alpha<1$, where $\mu$
and $\nu$ are defined by~\eqref{defmunu}.
\end{lemma}

We keep calling this inequality, ``Hardy-Poincar\'e inequality''
to remind the link with the inequality proved
in~\cite{BBDGV-cras,BBDGV-arch} but here we only use the
Poincar\'e type part of the inequality. The proof of the
Hardy-Poincar\'e inequality was performed in~\cite{BBDGV-cras} and
an estimate of the constant $\tilde\beta$ is also established. This
proof can be adapted to Lemma~\ref{lemmahardypoincare}. For
completeness, we give here the proof of this variant of the
Hardy-Poincar\'e inequality.\smallskip

\proof We first observe that we can reduce to show the inequality
for the Schwartz class $g\in{\cal D}(\RR^n)$ by simple
approximation arguments. In $L^2(\RR^n,\!\!\!\dd \mu)$ we consider
the closable quadratic form $v \mapsto \mathcal
Q(v):=\int_{\RR^n}|\nabla g|^2\dd \nu$ and $-\mathcal L$ the
unique non-negative, self-adjoint operator in $L^2(\RR^n,\!\!\dd
\mu)$ associated with the closure of $\mathcal Q$. By Persson's
theorem \cite{Pe}
\begin{equation*}
  \inf \sigma_{\rm ess}(- \mathcal L)= \lim_{R \to \infty} \inf_{v \in \mathcal H_R} \frac{\mathcal Q(v)}{\int_{\RR^n}|g|^2\dd \mu}
\end{equation*}
where $\mathcal H_R:=\{v \in H^1(\RR^d,\dd \nu)\,:\, {\rm supp}(v)
\subset \RR^n\setminus B(0,R)\}$. Roughly speaking it means that
the inequality is true for any weights with the same behaviour in
a neighbourhood of $+\infty$. By a straightforward computation
using~\eqref{equilibriumeq} and~\eqref{barenblattprofilefast}, we
have that
\[
\dd\mu(x)=\frac{1}{m}\left(D_*+\frac{1-\gamma}{mq}|x|^q\right)^{\frac{2-\gamma}{\gamma-1}}\dx
\sim_{|x|\to\infty}
\frac1m\left(\frac{1-\gamma}{m\,q}\right)^{\frac{2-\gamma}{\gamma-1}}|x|^{2\alpha-2}\dx
\]
and
\[
\dd\nu(x)=|x|^{2-q}\left(D_*+\frac{1-\gamma}{mq}|x|^q\right)^{\frac{1}{\gamma-1}}\dx
\sim_{|x|\to\infty}
\left(\frac{1-\gamma}{m\,q}\right)^{\frac{1}{\gamma-1}}|x|^{2\alpha}\dx
\]
with $\alpha$ chosen in such a way that
$q(2-\gamma)/(\gamma-1)=2(\alpha-1)$ and
$2-q+q/(\gamma-1)=2\alpha$, that is,
\begin{equation*}
\alpha=1+\frac{q(2-\gamma)}{2(\gamma-1)} \quad \mbox{or
equivalently} \quad \alpha=\frac{2-q}{2}+\frac{q}{2(\gamma-1)}.
\end{equation*}
It is left to the reader to check that $\alpha < -(d-2)/2$ in the
range of $1<p<\infty$ and $m_c<m$, and thus, we can apply
\cite[Theorem 1]{BBDGV-cras} to obtain
\begin{equation*}
  \inf \sigma_{\rm ess}(- \mathcal L) \ge \frac{1-\gamma}{q} \kappa_\alpha
\end{equation*}
where $\kappa_\alpha$ is the constant of the following Hardy
inequality with weight. We refer to \cite[Theorem 1]{BBDGV-cras}
for estimates on the constant $\kappa_\alpha$ depending if
$-n<\alpha$ or $\alpha\leq -n$. We remark that both cases happen
depending on the precise values of $m_c<m<\frac{n-p+1}{n(p-1)}$
and $1<p<\infty$.

The lowest eigenvalue of $-\mathcal L$ is $\lambda_1=0$ with
eigenfunctions given by the constants functions. Zero mean-value
solutions belong to the orthogonal set to the eigenspace
associated to $\lambda_1$. Since $\lambda_1$ is non-degenerate we
obtain the desired result for some $\lambda_2 \in
(0,\kappa_\alpha]$. \finproof

Since the only behaviour of the weights that counts is their growth
at infinity, we can avoid the singularity of the weight at the
origin for the singular case $1<p<2$ directly to obtain the
following stronger inequality.

\begin{corollary}
For any $\epsilon>0$, there exists a
constant $\tilde{\beta}_\epsilon>0$ such that
\begin{equation*}
\int_{\RR^n} g^2\dd\mu(x) \leq \tilde\beta_\epsilon\int_{\RR^n} |\nabla
g|^2\dd\nu_\epsilon(x),
\end{equation*}
for any function $g\in C^{1,\alpha}(\RR^n)$ satisfying
$\int_{\RR^n}g\dd\mu(x)=0$, where $\mu$ is defined
by~\eqref{defmunu}, and
\[ \dd\nu_\epsilon(x)=u_{D^*}(x)\left(\epsilon+|\nabla\left(F'\conv
u_{D^*}(x)\right)|\right)^{p-2}\dd x. \]

Therefore, setting $v=g/(F^{\prime\prime}\conv u_{D^*})$, we have the stronger linearised logarithmic Sobolev inequality
\begin{equation}\label{eqstronglogsob}
\EL\left[v\right] \leq \frac{\tilde\beta_\epsilon}{2}\IL^{\,(\epsilon)}\left[v\right]
\end{equation}
for all $v\in C^{1,\tilde\alpha}(\RR^n)$ with zero average and
$0<\tilde\alpha<1$, where
\begin{equation}\label{stronglinearfisherinfo}
I^{(\epsilon)} [v]= \int_{\RR^n} |\nabla\left(vF^{\prime\prime}\conv u_{D^*}\right)|^2\left(\epsilon+|\nabla\left(F'\conv
u_{D^*}(x)\right)|\right)^{p-2}u_{D^*}\dd x .
\end{equation}
\end{corollary}

\begin{remark}\label{bakry-emery}
The Bakry-Emery approach used in~{\rm \cite{CLMT}} to establish
the linearised logarithmic Sobolev
inequality~\eqref{linearizedlogsob} when $p=2$, does not seem to
apply here when $1<p<2$. For illustration, let us consider the
particular case $m=n=1$, that is the linearisation of the
1-dimensional rescaled $p$-Laplacian equation, $\partial_t v
=\mbox{\rm div}\,\{(p-1)u_{D^*}|A|^{p-2}W\}$. In this case, the
relative entropy dissipation
equation~\eqref{linearizeddissipation} simplifies as
\begin{equation*}
  \frac{\dd\EL\left[v(t)\right]}{\dd t} = -(p-1)\IL\left[v(t)\right],
\end{equation*}
and it is easy to show that its dissipation is
\[-\frac{\dd\IL\left[v(t)\right]}{\dd t} = 2(p-1)^2 \DD\left[v(t)|v_{D^*}\right],\]
where
\[\DD\left[v(t)|v_{D^*}\right]=\int_{\RR^n} F^{\prime\prime}\conv u_{D^*}\left[\mbox{\rm div}\left(u_{D^*}|A|^{p-2}W\right)\right]^2\dd x.\]

Following {\rm \cite{CLMT}}, if one can establish the estimate
$\DD\left[v(t)|v_{D^*}\right] \geq \lambda
\IL\left[v(t)|v_{D^*}\right]$ for some $\lambda>0$, then it will
imply the linearised logarithmic Sobolev
inequality~\eqref{linearizedlogsob}. But by a direct computation,
we can show that

\begin{eqnarray}\label{compareDI}
\DD\left[v(t)|v_{D^*}\right]  &=& \left(1+\frac{(p-2)^2}{p(p-1)}\right) \IL\left(v(t)|v_{D^*}\right) +
\int_{\RR^n} F^{\prime\prime}\circ u_{D^*}\left(u_{D^*} |A|^{p-2}\mbox{\rm div}\,W\right)^2\dd x \nonumber\\
& & + \; K\frac{p-2}{(p-1)^2}\int_{\RR^n} u_{D^*}|A|^{p-2}W^2|y|^{-q}\dd x.
\end{eqnarray}
Note that the second term in the above expression is non-negative
(because $F$ is convex), while the last term is non-positive in
the range $1<p< 2$. If $p\geq 2$, the last term is also
non-negative and we obtain $\DD\left(v(t)|v_{D^*}\right) \geq
\IL\left(v(t)|v_{D^*}\right)$, that is $\lambda=1$. This then
yields the desired inequality~\eqref{linearizedlogsob} when $p\geq
2$. But if $1<p<2$, one cannot derive this estimate
from~\eqref{compareDI}, at least at cursory glance.
\end{remark}

\noindent{\sl Proof of Theorem~{\rm \ref{theolinearizeddecay}\/ }}
We apply Lemma~\ref{lemmahardypoincare} to $g(t)=v(t)F^{\prime\prime}\!\circ\!u_{D^*}$ and $\beta=2(p-1)/\tilde\beta$. By  ~\eqref{linearizeddissipation} and ~\eqref{I} we have
\begin{equation*}
  \frac{\dd \EL\left[u(t)\right]}{\dd t} \leq -\beta\EL\left[v(t)\right].
\end{equation*}
This leads to~\eqref{linearizeddecay} by a Gronwall estimate.
\finproof
\section{Nonlinear stability}\label{sec:withrate}
\setcounter{equation}{0}

The first step to go from linear to nonlinear stability is to use
that our solution is sandwiched between two Barenblatt profiles to
compare the nonlinear relative entropy and its dissipation with
their linearised counterparts.

\begin{proposition}[Comparison linear/nonlinear relative entropy]\label{propcompareentropy}
Consider a function $u$ satisfying {\rm (H1)}. Then there exist
positive constants $C_1$ and $C_2$ such that
 \begin{equation*}
C_1\, \FL[u-u_{D_*}] \le \FF[u|u_{D_*}] \le C_2\,\FL[u-u_{D_*}]\;.
 \end{equation*}
\end{proposition}
\proof By Taylor's formula on the integrand of the
relative entropy we have
$$
F \conv u (x) - F \conv u_{D_*}(x) - F' \conv
u_{D_*}(x)\,(u(x)-u_{D_*}(x)) = \frac12 F''\!\circ\!\xi(x)
(u(x)-u_{D_*}(x))^2
$$
with
$$
u_{D_*}(x) W_0 \leq \min(u(x),u_{D_*}(x)) \leq \xi(x) \leq
\max(u(x),u_{D_*}(x)))\leq u_{D_*}(x) W_1,
$$
due to (H1), see also (H1'). The asserted result follows from homogeneity of $F''$ with $C_1:=mW_0^{\gamma-2}$ and $C_2:=mW_1^{\gamma-2}$.
\finproof

The next objective is to compare the nonlinear Fisher information,
$\II\left[u(t)|u_{D^*}\right]$, with its linear analogue,
$\IL\left[u(t)-u_{D^*}\right]$ along solutions of
\eqref{rescaledeq}. Let us point out that the weight
$$
|\nabla\left(F'\!\circ\!u_{D_*}\right)|^{p-2}=|x|^{2-q}
$$
in the linearised entropy dissipation diverges at the origin for
$1<p<2$. This singular behaviour makes complicated any attempt to
compare it with nonlinear Fisher analogues. Due to the singularity
of the weight $|\nabla \left(F^\prime\conv u_{D^*}\right)|^{p-2}$
at $x=0$, we will  replace $\IL\left[u(t)-u_{D^*}\right]$ by its
regularised analogue $\IL^{(\epsilon)}\left[u(t)-u_{D^*}\right]$
defined by (\ref{stronglinearfisherinfo}), where $\epsilon>0$ is a
fixed constant.

\begin{proposition}[Comparison linear/non-linear Fisher information] \label{propcomparefisher}
Assume that $u$ is the solution of~\eqref{rescaledeq}, and set
$v=u-u_{D^*}$. Then:
\begin{enumerate}
\item Case $1 < p < 2$: Given $\epsilon >0$, there exist $t_0>0$
and positive constants $\kappa_1$ and $\kappa_2$ such that for all
$t>t_0$,
\begin{equation}\label{comparefisher}
\IL^{\,(\epsilon)}\left[v(t)\right] \leq \kappa_1\,
\II\left[u(t)|u_{D^*}\right] + \kappa_2\FL\left[v(t)\right].
\end{equation}

\item Case $p>2$: There exist $t_0>0$ and positive constants
$\kappa_1$ and $\kappa_2$ such that for all $t>t_0$,
\begin{equation}\label{comparefisher2}
\IL\left[v(t)\right] \leq \kappa_1\, \II\left[u(t)|u_{D^*}\right]
+ \kappa_2\FL\left[v(t)\right].
\end{equation}
\end{enumerate}
Moreover $\kappa_2$ can be chosen arbitrary small provided that
$t_0$ is large enough.
\end{proposition}

\noindent The proof of this proposition is organised as follows:

\noindent{\bf Claim 1:} We first show that, for all $\epsilon \geq
0$, there exists $\kappa_0>0$ such that
\begin{equation}\label{eqstep1}
\IL^{\,(\epsilon)}\left[v(t)\right] \leq \kappa_0\IL_\gamma^{\,(\epsilon)}\left[v(t)\right] + \kappa_2\FL\left[v(t)\right],
\end{equation}
where
\begin{equation}\label{defIgammaepsilon}
\IL_\gamma^{\,(\epsilon)}[v] =\int_{\RR^n}
\left|\nabla\left[F'\conv u - F'\conv u_{D^*}\right]\right|^2
\left(\epsilon+|\nabla\left(F'\conv
u_{D^*}(x)\right)|\right)^{p-2}u_{D^*}\dd x .
\end{equation}

\noindent{\bf Claim 2:} Next we show that if $1<p<2$, then  for all
$\epsilon >0$ there exists $\delta>0$ such that
\begin{equation}\label{eqstep2}
\IL_\gamma^{\,(\epsilon)}\left[v(t)\right]  \leq \delta
\II\left[u(t)|u_{D^*}\right],
\end{equation}
whereas if $2<p<\infty$, then there exists $\delta>0$ such that
\begin{equation*}
\IL_\gamma\left[v(t)\right]  \leq \delta
\II\left[u(t)|u_{D^*}\right].
\end{equation*}
Combining (\ref{eqstep1}) and (\ref{eqstep2}), we obtain the
desired inequalities \eqref{comparefisher}-\eqref{comparefisher2}
with $\kappa_1=\delta\kappa_0$.

\

\noindent{\sl Proof of Claim 1\/:} Here we follow the arguments of
the proof of Lemma 5.1 in \cite{BBDGV-arch}. Indeed, let
$h_k(w)=(w^{k-1}-1)/(k-1)$, where
\begin{equation*}
  w(t,x)=\frac{u(t,x)}{u_{D^*}(x)}.
\end{equation*}
Because of assumption {\rm (H1)}, we have that $W_0\leq w(t,x)\leq
W_1$, where the constant $W_0$ and $W_1$ are such that $0<W_0 <
1 < W_1$. By studying the function $h_2/h_\gamma$ on $[W_0, W_1]$, we have
\begin{equation}\label{eqstep1-1}
\alpha_0 h_\gamma(w)^2 \leq h_2(w)^2 \leq \alpha_1 h_\gamma(w)^2,
\end{equation}
and
\begin{equation*}
 h'_2(w)^2 \leq \alpha_2 h'_\gamma(w)^2,
\end{equation*}
where
\[\alpha_0 :=|\gamma-1|^2\left|\frac{W_0-1}{W_0^{\gamma-1}-1}\right|^2 <1, \quad \alpha_1 := |\gamma-1|^2\left|\frac{W_1-1}{W_1^{\gamma-1}-1}\right|^2 >1\]
and
\[\alpha_2 := W_1^{2(2-\gamma)} > 1.\]
Now, define
\[
\IL_k^{\,(\epsilon)}[v] := m^2\int_{\RR^n}
\left|\nabla\left(u_{D^*}^{\gamma-1}h_k(w)\right)\right|^2(\epsilon+|x|^{q-1})^{p-2}u_{D^*}
\dd x.
\]
We have that $ \IL^{\,(\epsilon)}=\IL_2^{\,(\epsilon)}$ and for
$k=\gamma$,  $\IL_\gamma^{\,(\epsilon)}$ is defined in~\eqref{defIgammaepsilon}. Next we compute
$\IL_k^{\,(\epsilon)}[v]$. By expanding
$\left|\nabla\left(u_{D^*}^{\gamma-1}h_k(w)\right)\right|^2$, we
have
\begin{eqnarray*}
\IL_k^{\,(\epsilon)}[v]  &=& m^2\int_{\RR^n}u_{D^*}^{2\gamma-1}h'_k(w)^2|\nabla w|^2(\epsilon+|x|^{q-1})^{p-2}\dd x \\
& & + (1-\gamma)^2\int_{\RR^n} h_k(w)^2|x|^{2(q-1)}(\epsilon+|x|^{q-1})^{p-2}u_{D^*}\dd x\\
& & +2m(1-\gamma)\int_{\RR^n} u_{D^*}^\gamma h'_k(w)h_k(w)|x|^{q-2}(\epsilon+|x|^{q-1})^{p-2}\nabla w\cdot x\dd x.
\end{eqnarray*}
Integrating by parts, the last integral can be rewritten as
\begin{align*}
\int_{\RR^n} u_{D^*}^\gamma
h'_k(w)h_k(w)|x|^{q-2}(\epsilon+|x|^{q-1})^{p-2}&\nabla w\cdot
x\dd x \\
=&\, \frac{1}{2}\int_{\RR^n} \nabla\left(h_k(w)^2\right) \cdot |x|^{q-2}x(\epsilon+|x|^{q-1})^{p-2} u_{D^*}^\gamma\dd x\\
=&\, -\frac{1}{2}\int_{\RR^n} h_k(w)^2\mbox{div}\left(|x|^{q-2}x(\epsilon+|x|^{q-1})^{p-2} u_{D^*}^\gamma\right)\dd x\\
=&\, -\frac{1}{2}\int_{\RR^n} h_k(w)^2\mbox{div}\left(|x|^{q-2}x(\epsilon+|x|^{q-1})^{p-2}\right)u_{D^*}^\gamma\dd x\\
 &\, +\frac{\gamma}{2m}\int_{\RR^n} h_k(w)^2|x|^{2(q-1)}(\epsilon+|x|^{q-1})^{p-2} u_{D^*}\dd x.
\end{align*}
Then,
\begin{align}\label{valueIkepsilon}
\IL_k^{\,(\epsilon)}[v] =&\, m^2\int_{\RR^n}u_{D^*}^{2\gamma-1}h'_k(w)^2|\nabla w|^2(\epsilon+|x|^{q-1})^{p-2}\dd x \nonumber\\
& + (1-\gamma)\int_{\RR^n}
h_k(w)^2|x|^{2(q-1)}(\epsilon+|x|^{q-1})^{p-2}u_{D^*}\dd x
\nonumber\\ &
-m(1-\gamma)\int_{\RR^n}h_k(w)^2\mbox{div}\left(|x|^{q-2}x(\epsilon+|x|^{q-1})^{p-2}\right)u_{D^*}^\gamma
\dd x.
\end{align}
Next we set $\kappa_0:=\max(\alpha_1, \alpha_2)$. Moreover, since $w$ uniformly converges to 1 as $t$ goes to $\infty$, then $\alpha_0$, $\alpha_1$,
$\alpha_2$ and $\kappa_0>1$ can be chosen arbitrary close to $1$ provided that
$t>t_0$, for some $t_0$ large enough. Combining \eqref{eqstep1-1}-\eqref{valueIkepsilon}, we have
\begin{align*}
\IL^{\,(\epsilon)}[v] =\IL_2^{\,(\epsilon)}[v] \leq&\,
m^2\alpha_2\int_{\RR^n}u_{D^*}^{2\gamma-1}h'_\gamma(w)^2 |\nabla
w|^2(\epsilon+|x|^{q-1})^{p-2}\dd x
\\
& + \alpha_1(1-\gamma)\int_{\RR^n} h_\gamma(w)^2|x|^{2(q-1)}(\epsilon+|x|^{q-1})^{p-2}u_{D^*}\dd x \\
& -m(1-\gamma)\int_{\RR^n}h_2(w)^2\mbox{div}\left(|x|^{q-2}x(\epsilon+|x|^{q-1})^{p-2}\right)u_{D^*}^\gamma\dd x\\
&\hspace{-1.7cm}\leq\, \kappa_0 \IL_\gamma^{\,(\epsilon)}[v] +
m(1-\gamma) \int_{\RR^n} u_{D^*}^\gamma
\mbox{div}\left(|x|^{q-2}x(\epsilon+|x|^{q-1})^{p-2}\right)\left(\kappa_0
h_\gamma(w)^2-h_2(w)^2\right) \dd x.
\end{align*}
Finally, we observe that $0\leq \kappa_0 h_\gamma(w)^2-h_2(w)^2 \leq \left(\kappa_0/\alpha_0-1\right)h_2(w)^2$ and by direct computation
\begin{equation}\label{estimatediv}
\left|\mbox{ div}\left(|x|^{q-2}x(\epsilon+|x|^{q-1})^{p-2}\right)\right|
\leq n+2(q-2),
\end{equation}
and we then deduce that
\[\IL^{\,(\epsilon)}[v(t)] \leq \kappa_0 \IL_\gamma^{\,(\epsilon)}[v(t)] + \kappa_2\FL\left[v(t)\right] \]
with $\kappa_2:=2(\kappa_0/\alpha_0-1)(1-\gamma)(n+2(q-2))>0$. Clearly $\kappa_2$ is arbitrary small provided that $t_0$ is large enough.
\finproof

\noindent{\sl Proof of Claim 2\/; case $1<p<2$:} First we expand
$\IL_\gamma^{\,(\epsilon)}[v(t)]$ and
$\II\left[u(t)|u_{D^*}\right]$, and we have that
\begin{align*}
\IL_\gamma^{\,(\epsilon)}[v(t)]  =&\,   \int_{\RR^n} |\nabla\left[F'\circ u\right]|^2 (\epsilon+|\nabla\left[F'\circ u_{D^*}\right]|)^{p-2}u_{D^*}\dd x  \\
& +  \int_{\RR^n}|\nabla\left[F'\circ u_{D^*}\right]|^2(\epsilon+|\nabla\left[F'\circ u_{D^*}\right]|)^{p-2}u_{D^*}\dd x\\
& - 2\int_{\RR^n}\nabla\left[F'\circ u\right] \cdot
\nabla\left[F'\circ u_{D^*}\right](\epsilon+|\nabla\left[F'\circ
u_{D^*}\right]|)^{p-2}u_{D^*}\dd x
\end{align*}
and
\begin{align*}
\II\left[u(t)|u_{D^*}\right] =&\,  \int_{\RR^n}
\left|\nabla\left[F'\circ u\right]\right|^pu\dd x  +
\int_{\RR^n}|\nabla\left[F'\circ u_{D^*}\right]|^pu\dd x \\ & -
\int_{\RR^n}\nabla\left[F'\circ u\right] \cdot \nabla
c^*\left(\nabla\left[F'\circ u_{D^*}\right] \right)u\dd x \\ &  -
\int_{\RR^n}\nabla\left[F'\circ u_{D^*}\right] \cdot \nabla
c^*\left(\nabla\left[F'\circ u\right]\right)u\dd x.
\end{align*}
Next we use Young inequality $a\cdot b \leq c(a) + c^*(b)$ with
$c(z)=|z|^q/q$, $a= \nabla c^*\left(\nabla\left[F'\circ
u\right]\right)$ and $b=\nabla\left[F'\circ u_{D^*}\right]$, to
have that
\[
\nabla\left[F'\circ u_{D^*}\right] \cdot \nabla c^*\left(\nabla
F'(u)\right) \leq \frac{1}{q} \left|\nabla\left[F'\circ
u\right]\right|^p + \frac{1}{p} |\nabla\left[F'\circ
u_{D^*}\right]|^p.
\]
Then $\II\left[u(t)|u_{D^*}\right]$ can be estimated as
\begin{align}\label{estimate1}
\II\left[u(t)|u_{D^*}\right] \geq&\, \frac{1}{p}
\int_{\RR^n}\left|\nabla\left[F'\circ u\right]\right|^pu \dd x +
\frac{1}{q}
\int_{\RR^n} \left|\left(\nabla\left[F'\circ u_{D^*}\right]\right)\right|^p u \dd x\nonumber \\
&- \int_{\RR^n} \nabla\left[F'\circ u\right] \cdot \nabla c^*\left(\nabla\left[F'\circ u_{D^*}\right]\right)u \dd x \nonumber\\
=&\, \frac{1}{p} \int_{\RR^n}\left|\nabla\left[F'\circ
u\right]\right|^pu\dd x + \frac{1}{q}
\int_{\RR^n}|\left(\nabla\left[F'\circ
u_{D^*}\right]\right)|^pu\dd x \nonumber \\
&- - \!\frac{nm}{\gamma(1-\gamma)}\int_{\RR^n} u^\gamma\dd x.
\end{align}
Now, we compute the cross term of $\IL_\gamma^{\,(\epsilon)}[v(t)]$. We have that
\begin{align*}
\int_{\RR^n}\nabla\left[F'\circ u\right]) \cdot \nabla\left[F'\circ u_{D^*}\right](\epsilon+&|\nabla\left[F'\circ u_{D^*}\right]|)^{p-2}u_{D^*}\dd x\\
=&\, \frac{m}{1-\gamma} \int_{\RR^n} \nabla(u^{\gamma-1})\cdot x|x|^{q-2}(\epsilon+|x|^{q-1})^{p-2}u_{D^*}\dd x\\
=&\, -\frac{m}{1-\gamma} \int_{\RR^n} u^{\gamma-1} \, \mbox{div}\left(x|x|^{q-2}(\epsilon+|x|^{q-1})^{p-2}u_{D^*}\right)\dd x\\
=&\, -\frac{m}{1-\gamma}\int_{\RR^n} \frac{u^\gamma}{w} \, \mbox{div}\left(x|x|^{q-2}(\epsilon+|x|^{q-1})^{p-2}\right)\dd x\\
& +\frac{1}{1-\gamma} \int_{\RR^n}
w^{\gamma-1}u_{D^*}|x|^{2(q-1)}(\epsilon+|x|^{q-1})^{p-2}\dd x.
\end{align*}
Since the last term in the above sum is non-negative, then using
(\ref{estimatediv}) and the fact that $1<p< 2$, we can estimate
$\IL_\gamma^{\,(\epsilon)}[v(t)] $ as
\begin{eqnarray*}
\IL_\gamma^{\,(\epsilon)}[v(t)]  &\leq&  \int_{\RR^n} \left|\nabla\left[F'\circ u\right]\right|^2 (\epsilon+|\nabla\left[F'\circ u_{D^*}\right]|)^{p-2}u_{D^*}\dd x \\
& &+  \int_{\RR^n}|\nabla\left[F'\circ u_{D^*}\right]|^2(\epsilon+|\nabla\left[F'\circ u_{D^*}\right]|)^{p-2}u_{D^*}\dd x\\
& & +\frac{2m}{1-\gamma} \int_{\RR^n} \frac{u^\gamma}{w}\mbox{div}\left(x|x|^{q-2}(\epsilon+|x|^{q-1})^{p-2}\right)\dd x  \\
&\leq& \int_{\RR^n} \left|\nabla\left[F'\circ u\right]\right|^2 (\epsilon+|\nabla\left[F'\circ u_{D^*}\right]|)^{p-2}u_{D^*}\dd x  \\
& & +  \int_{\RR^n}|\nabla\left[F'\circ u_{D^*}\right]|^2(\epsilon+|\nabla\left[F'\circ u_{D^*}\right]|)^{p-2}u_{D^*}\dd x\\
& & +\frac{2m(n+2(q-2))}{1-\gamma} \int_{\RR^n} \frac{u^\gamma}{w} \dd x\nonumber\\
& \leq & \int_{\RR^n} \left|\nabla\left[F'\circ u\right]\right|^p
\Phi_\epsilon(u,u_{D^*})^{2-p} u_{D^*}\dd x
+ \int_{\RR^n}|\nabla\left[F'\circ u_{D^*}\right]|^p u_{D^*}\dd x \\
& & + \frac{2m(n+2(q-2))}{1-\gamma} \int_{\RR^n} \frac{u^\gamma}{w} \dd x,
\end{eqnarray*}
where
\[\Phi_\epsilon(u, u_{D^*}) := \frac{\left|\nabla\left[F'\circ u\right]\right|}{\epsilon+|\nabla\left[F'\circ u_{D^*}\right]|}. \]
Using Lemma \ref{ckesti} and \eqref{newbehinfty}, the reader can
check that $\Phi_\epsilon(u, u_{D^*})$ is uniformly bounded by
expressing the gradients in terms of derivatives of $w$ and
$u_{D^*}$, and thus, $\Phi_\epsilon(u, u_{D^*}) \leq \eta$ for
some $\eta>0$ depending on $\epsilon$. Therefore, we obtain the
estimate
\begin{eqnarray}\label{estimate2}
\IL_\gamma^{\,(\epsilon)}[v(t)]  &\leq&  \frac{\eta^{2-p}}{W_0} \int_{\RR^n}
\left|\nabla\left[F'\circ u\right]\right|^p u\dd x + \frac{1}{W_0}\int_{\RR^n} \left|\nabla\left[F'\circ u_{D^*}\right]\right|^p u\dd x \nonumber \\
& & + \frac{2m(n+2(q-2))}{(1-\gamma)W_0} \int_{\RR^n} u^\gamma \dd x.
\end{eqnarray}
Combining~\eqref{estimate1} and~\eqref{estimate2}, and setting
\[\delta := \frac{1}{W_0}\max\left(p\eta^{2-p}, q, \frac{2\gamma(n+2(q-2))}{n} \right),\]
we deduce~\eqref{eqstep2}.
\finproof

\noindent{\sl Proof of Claim 2\/; case $2<p<\infty$:} As above, we
have the expression
$$
\IL_\gamma\left[v(t)\right] = \int_{\RR^n}
\left|\nabla\left[F'\conv u - F'\conv u_{D^*}\right]\right|^2
|\nabla\left(F'\conv u_{D^*}(x)\right)|^{p-2}u_{D^*}\dd x .
$$
For convenience, we can also rewrite
$\II\left(u(t)|u_{D^*}\right)$ as
\begin{equation*}
\II\left[u(t)|u_{D^*}\right]=\int_{\cal K}
H[u|u_{D^*}]\,\left|\nabla\left[F'\conv u-F'\conv
u_{D^*}\right]\right|^2 u \dd y,
\end{equation*}
where
\begin{equation*}
H[u(t)|u_{D^*}]=\frac{\nabla\left(F'\conv u-F'\conv
u_{D^*}\right)\cdot\left[\nabla c^*\left[\nabla \left(F'\conv
u\right)\right]-\nabla c^*\left[\nabla \left(F'\conv
u_{D^*}\right)\right]\right]}{|\nabla\left(F'\conv u-F'\conv
u_{D^*}\right)|^2}
\end{equation*}
and ${\cal K}:=\{x\in\RR^n \mbox{ such that } |\nabla\left[F'\conv
u-F'\conv u_{D^*}\right]|\neq 0 \}$. Let us show that there exist
a constant $\delta>0$, such that for all $t>t_0$,
\begin{equation}\label{chogella}
H[u(t)|u_{D^*}]\geq \delta\,|\nabla\left(F'\conv
u_{D^*}\right)|^{p-2}.
\end{equation}
Let us remark, if $p=2$, then $\delta=1$, and equality holds
in~\eqref{chogella}.

For simplicity, set $a(t)=\nabla\left(F'\conv u(t)\right)$ and
$a_{D^*}=\nabla\left(F'\conv u_{D^*}\right)$. It is clear that
\eqref{chogella} holds in the set where $a_{D^*}=0$. Therefore,
let us restrict to the set where $a_{D^*}\neq 0$ without loss of
generality. Let us denote $b(t)=a(t)/|a_{D^*}|$ and
$b_{D^*}=a_{D^*}/|a_{D^*}|$. It is straightforward to check that
\begin{equation}\label{eqn14}
\frac{H[u(t)|u_{D^*}]}{|\nabla\left(F'\conv u_{D^*}\right)|^{p-2}}
=
\frac{\left(b(t)-b_{D^*}\right)\cdot\left(|b(t)|^{p-2}b(t)-b_{D^*}\right)}{|b(t)-b_{D^*}|^2}.
\end{equation}
Let $\theta$ denote the angle between $b(t)$ and $b_{D^*}$. We
have that
\[|b-b_{D^*}|^2 =|b|^2+|b_{D^*}|^2-2b\cos\theta = 1+|b|^2-2b\cos\theta,\]
and
\[(b-b_{D^*})\cdot \left(|b|^{p-2}b-b_{D^*})\right)=|b|^p-|b|\cos\theta-|b|^{p-1}\cos\theta+|b_{D^*}|^2=1+|b|^p-(|b|+|b|^{p-1})\cos\theta\]
so that~\eqref{eqn14} reads as:
\begin{equation*}
\frac{H[u(t)|u_{D^*}]}{|\nabla\left(F'\conv u_{D^*}\right)|^{p-2}}
=\frac{1+r(t)^p-\left(r(t)+r(t)^{p-1}\right)x(t)}{1+r(t)^2-2r(t)x(t)}
\end{equation*}
where $r(t)=|b(t)|\geq 0$ and $x(t)=\cos\theta\in [-1,1]$, with
$r(t)\to 1$ as $t\to \infty$. Estimate \eqref{chogella} is reduced
to show that
\begin{equation}\label{eqn16}
f_p(r,x) :=\frac{1+r^p-(r+r^{p-1})x}{1+r^2-2rx} \geq \delta,
\end{equation}
for all $x\in[-1,1]$ and for all $r\geq 0$. For that, let us
define the function
$$
F_p(r,x)=:= 1+r^p-(r+r^{p-1})x - \delta (1+r^2-2rx),
$$
which is easily checked to be decreasing in $x$ for $r>0$,
whenever $\delta\leq \frac12$. Therefore, we have $F_p(r,x)>
F_p(r,1)$ and thus, to show \eqref{eqn16} for $r>0$ and $-1\leq
x<1$ is reduced to show that $F_p(r,1)\geq 0$, whenever
$\delta\leq \frac12$. Since $F_p(1,1)=0$, this is equivalent to
show that $f_p(r,1)\geq \delta$ for $0<r<1$ and $r>1$. The last
assertion comes from the fact that when $p>2$,
\[
f_p(r,1)=\frac{(r^p-r)-(r^{p-1}-1)}{(r-1)^2}= \frac{(r-1)(r^{p-1}-1)}{(r-1)^2}=\frac{r^{p-1}-1}{r-1}
\]
is bounded below by 1, since in $0<r<1$, we have $r^{p-1}-1<r-1<0$; in $r>1$, we have 
$r^{p-1}-1>r-1>0$;  and $\lim_{r\to 1} f_p(r,1)=p-1>1$. Therefore,
\[f_p(r,1) \geq 1 > \frac{1}{2} \geq \delta.\]
\finproof

\noindent{\sl Proof of the main theorem,
Theorem~\ref{maintheo}\/:} Given $1<p<2$ and $\epsilon> 0$, set
$v=u-u_{D^*}$. From Proposition~\ref{propcomparefisher} and the
strong linearised logarithmic Sobolev
inequality~\eqref{eqstronglogsob}, we have that
\[\EL\left[v(t)\right] \leq \frac{\kappa_1\tilde{\beta}_\epsilon}{2-\kappa_2\tilde{\beta}_\epsilon} \,\II\left[u(t)|u_{D^*}\right],\]
where $\kappa_2$ can be chosen arbitrary small provided that $t>t_0$ is large enough. This together with proposition \ref{propcompareentropy} yields the logarithmic Sobolev
type inequality:
\begin{equation}\label{nonlinearlogsob}
\EE\left[u(t)|u_{D^*}\right] \leq \frac{1}{\lambda}
\,\II\left[u(t)|u_{D^*}\right],
\end{equation}
where $\lambda:=\left(2-\kappa_2\tilde{\beta}_\epsilon\right)/C_2\kappa_1\tilde{\beta}_\epsilon>0$.
We combine (\ref{nonlinearlogsob}) and the entropy dissipation equation
\[  \frac{\dd}{\dd t}  \FF[u(t)|u_{D_*}] =- \II[u(t)|u_{D_*}] \]
to obtain the exponential decay of the relative entropy, $\EE\left[u(t)|u_{D^*}\right] \leq e^{-\lambda t} \EE\left[u_0|u_{D^*}\right] $. The $L^1$-decay (\ref{recaledL1decay}) follows from the Csisz\`ar-Kullback type inequality (see for e.g., \cite{Ag-na08}),
\[ \|u(t)-u_{D^*}\|^2_{L^1(\RR^n)} \leq M(n,n,p) \EE\left[u(t)|u_{D^*}\right], \;\; M(m,n,p)>0, \]
 and (\ref{L1decay}) is a direct consequence of the rescaling
 (\ref{rescaledrho})-(\ref{defR}). The case $p>2$ follows
 analogously without need of using the regularised entropy dissipation.
\finproof
\begin{remark}
  The rate of convergence $\lambda$ in~Theorem~\ref{maintheo} can be explicitely reconstructed in the above computation for a given choice of $\epsilon$ for $p<2$, but $t_0$ will depend on the choice of $\epsilon$. For $p>2$, we can even give an explicit range for the constant $\lambda$.
\end{remark}

\noindent {\bf Acknowledgements:} MA is partially supported by a
grant from NSERC (Canada). JAC acknowledge partial support from
the project MTM2008-06349-C03-03 from DGI-MCI (Spain). We also
thank WPI-Vienna where this research was started in the framework
of the program "Optimal transportation structures, gradient flows
and entropy methods for Applied PDEs". AB and JAC acknowledges
IPAM--UCLA where part of this work was done.

\bigskip\noindent{\small This paper is under the Creative Commons licence Attribution-NonCommercialShareAlike~2.5.}

\end{document}